\documentclass[11pt,twoside,a4paper]{article}
\usepackage{makeidx}
\usepackage[dvips]{graphicx}
\usepackage{amscd}
\usepackage{amsmath}
\usepackage{amssymb}
\usepackage{latexsym}
\usepackage[latin1]{inputenc}
\usepackage[T1]{fontenc}   
\usepackage[english]{babel}
\newcommand{\tun}{\begin{picture}(5,0)(-2,-1)
\put(0,0){\circle*{2}}
\end{picture}}

\newcommand{\tdeux}{\begin{picture}(7,7)(0,-1)
\put(3,0){\circle*{2}}
\put(3,0){\line(0,1){5}}
\put(3,5){\circle*{2}}
\end{picture}}

\newcommand{\ttroisun}{\begin{picture}(15,8)(-5,-1)
\put(3,0){\circle*{2}}
\put(-0.65,0){$\vee$}
\put(6,7){\circle*{2}}
\put(0,7){\circle*{2}}
\end{picture}}
\newcommand{\ttroisdeux}{\begin{picture}(5,12)(-2,-1)
\put(0,0){\circle*{2}}
\put(0,0){\line(0,1){5}}
\put(0,5){\circle*{2}}
\put(0,5){\line(0,1){5}}
\put(0,10){\circle*{2}}
\end{picture}}

\newcommand{\tquatreun}{\begin{picture}(15,12)(-5,-1)
\put(3,0){\circle*{2}}
\put(-0.65,0){$\vee$}
\put(6,7){\circle*{2}}
\put(0,7){\circle*{2}}
\put(3,7){\circle*{2}}
\put(3,0){\line(0,1){7}}
\end{picture}}
\newcommand{\tquatredeux}{\begin{picture}(15,18)(-5,-1)
\put(3,0){\circle*{2}}
\put(-0.65,0){$\vee$}
\put(6,7){\circle*{2}}
\put(0,7){\circle*{2}}
\put(0,14){\circle*{2}}
\put(0,7){\line(0,1){7}}
\end{picture}}
\newcommand{\tquatretrois}{\begin{picture}(15,18)(-5,-1)
\put(3,0){\circle*{2}}
\put(-0.65,0){$\vee$}
\put(6,7){\circle*{2}}
\put(0,7){\circle*{2}}
\put(6,14){\circle*{2}}
\put(6,7){\line(0,1){7}}
\end{picture}}
\newcommand{\tquatrequatre}{\begin{picture}(15,18)(-5,-1)
\put(3,5){\circle*{2}}
\put(-0.65,5){$\vee$}
\put(6,12){\circle*{2}}
\put(0,12){\circle*{2}}
\put(3,0){\circle*{2}}
\put(3,0){\line(0,1){5}}
\end{picture}}
\newcommand{\tquatrecinq}{\begin{picture}(9,19)(-2,-1)
\put(0,0){\circle*{2}}
\put(0,0){\line(0,1){5}}
\put(0,5){\circle*{2}}
\put(0,5){\line(0,1){5}}
\put(0,10){\circle*{2}}
\put(0,10){\line(0,1){5}}
\put(0,15){\circle*{2}}
\end{picture}}

\newcommand{\tcinqun}{\begin{picture}(20,8)(-5,-1)
\put(3,0){\circle*{2}}
\put(-0.5,0){$\vee$}
\put(6,7){\circle*{2}}
\put(0,7){\circle*{2}}
\put(3,0){\line(2,1){10}}
\put(3,0){\line(-2,1){10}}
\put(-7,5){\circle*{2}}
\put(13,5){\circle*{2}}
\end{picture}}
\newcommand{\tcinqdeux}{\begin{picture}(15,14)(-5,-1)
\put(3,0){\circle*{2}}
\put(-0.65,0){$\vee$}
\put(6,7){\circle*{2}}
\put(0,7){\circle*{2}}
\put(3,7){\circle*{2}}
\put(3,0){\line(0,1){7}}
\put(0,7){\line(0,1){7}}
\put(0,14){\circle*{2}}
\end{picture}}
\newcommand{\tcinqtrois}{\begin{picture}(15,15)(-5,-1)
\put(3,0){\circle*{2}}
\put(-0.65,0){$\vee$}
\put(6,7){\circle*{2}}
\put(0,7){\circle*{2}}
\put(3,7){\circle*{2}}
\put(3,0){\line(0,1){7}}
\put(3,7){\line(0,1){7}}
\put(3,14){\circle*{2}}
\end{picture}}
\newcommand{\tcinqquatre}{\begin{picture}(15,14)(-5,-1)
\put(3,0){\circle*{2}}
\put(-0.65,0){$\vee$}
\put(6,7){\circle*{2}}
\put(0,7){\circle*{2}}
\put(3,7){\circle*{2}}
\put(3,0){\line(0,1){7}}
\put(6,7){\line(0,1){7}}
\put(6,14){\circle*{2}}
\end{picture}}
\newcommand{\tcinqcinq}{\begin{picture}(15,19)(-5,-1)
\put(3,0){\circle*{2}}
\put(-0.65,0){$\vee$}
\put(6,7){\circle*{2}}
\put(0,7){\circle*{2}}
\put(6,14){\circle*{2}}
\put(6,7){\line(0,1){7}}
\put(0,14){\circle*{2}}
\put(0,7){\line(0,1){7}}
\end{picture}}
\newcommand{\tcinqsix}{\begin{picture}(15,20)(-7,-1)
\put(3,0){\circle*{2}}
\put(-0.65,0){$\vee$}
\put(6,7){\circle*{2}}
\put(0,7){\circle*{2}}
\put(-3.65,7){$\vee$}
\put(3,14){\circle*{2}}
\put(-3,14){\circle*{2}}
\end{picture}}
\newcommand{\tcinqsept}{\begin{picture}(15,8)(-5,-1)
\put(3,0){\circle*{2}}
\put(-0.65,0){$\vee$}
\put(6,7){\circle*{2}}
\put(0,7){\circle*{2}}
\put(2.35,7){$\vee$}
\put(3,14){\circle*{2}}
\put(9,14){\circle*{2}}
\end{picture}}
\newcommand{\tcinqhuit}{\begin{picture}(15,26)(-5,-1)
\put(3,0){\circle*{2}}
\put(-0.65,0){$\vee$}
\put(6,7){\circle*{2}}
\put(0,7){\circle*{2}}
\put(0,14){\circle*{2}}
\put(0,7){\line(0,1){7}}
\put(0,21){\circle*{2}}
\put(0,14){\line(0,1){7}}
\end{picture}}
\newcommand{\tcinqneuf}{\begin{picture}(15,26)(-5,-1)
\put(3,0){\circle*{2}}
\put(-0.65,0){$\vee$}
\put(6,7){\circle*{2}}
\put(0,7){\circle*{2}}
\put(6,14){\circle*{2}}
\put(6,7){\line(0,1){7}}
\put(6,21){\circle*{2}}
\put(6,14){\line(0,1){7}}
\end{picture}}
\newcommand{\tcinqdix}{\begin{picture}(15,19)(-5,-1)
\put(3,5){\circle*{2}}
\put(-0.5,5){$\vee$}
\put(6,12){\circle*{2}}
\put(0,12){\circle*{2}}
\put(3,0){\circle*{2}}
\put(3,0){\line(0,1){12}}
\put(3,12){\circle*{2}}
\end{picture}}
\newcommand{\tcinqonze}{\begin{picture}(15,26)(-5,-1)
\put(3,5){\circle*{2}}
\put(-0.65,5){$\vee$}
\put(6,12){\circle*{2}}
\put(0,12){\circle*{2}}
\put(3,0){\circle*{2}}
\put(3,0){\line(0,1){5}}
\put(0,12){\line(0,1){7}}
\put(0,19){\circle*{2}}
\end{picture}}
\newcommand{\tcinqdouze}{\begin{picture}(15,26)(-5,-1)
\put(3,5){\circle*{2}}
\put(-0.65,5){$\vee$}
\put(6,12){\circle*{2}}
\put(0,12){\circle*{2}}
\put(3,0){\circle*{2}}
\put(3,0){\line(0,1){5}}
\put(6,12){\line(0,1){7}}
\put(6,19){\circle*{2}}
\end{picture}}
\newcommand{\tcinqtreize}{\begin{picture}(5,26)(-2,-1)
\put(0,0){\circle*{2}}
\put(0,0){\line(0,1){7}}
\put(0,7){\circle*{2}}
\put(0,7){\line(0,1){7}}
\put(0,14){\circle*{2}}
\put(-3.65,14){$\vee$}
\put(-3,21){\circle*{2}}
\put(3,21){\circle*{2}}
\end{picture}}
\newcommand{\tcinqquatorze}{\begin{picture}(9,26)(-5,-1)
\put(0,0){\circle*{2}}
\put(0,0){\line(0,1){5}}
\put(0,5){\circle*{2}}
\put(0,5){\line(0,1){5}}
\put(0,10){\circle*{2}}
\put(0,10){\line(0,1){5}}
\put(0,15){\circle*{2}}
\put(0,15){\line(0,1){5}}
\put(0,20){\circle*{2}}
\end{picture}}

\newcommand{\tdun}[1]{\begin{picture}(10,5)(-2,-1)
\put(0,0){\circle*{2}}
\put(3,-2){\tiny #1}
\end{picture}}

\newcommand{\tddeux}[2]{\begin{picture}(12,5)(0,-1)
\put(3,0){\circle*{2}}
\put(3,0){\line(0,1){5}}
\put(3,5){\circle*{2}}
\put(6,-2){\tiny #1}
\put(6,3){\tiny #2}
\end{picture}}

\newcommand{\tdtroisun}[3]{\begin{picture}(20,12)(-5,-1)
\put(3,0){\circle*{2}}
\put(-0.65,0){$\vee$}
\put(6,7){\circle*{2}}
\put(0,7){\circle*{2}}
\put(5,-2){\tiny #1}
\put(9,5){\tiny #2}
\put(-5,5){\tiny #3}
\end{picture}}
\newcommand{\tdtroisdeux}[3]{\begin{picture}(12,12)(-2,-1)
\put(0,0){\circle*{2}}
\put(0,0){\line(0,1){5}}
\put(0,5){\circle*{2}}
\put(0,5){\line(0,1){5}}
\put(0,10){\circle*{2}}
\put(3,-2){\tiny #1}
\put(3,3){\tiny #2}
\put(3,9){\tiny #3}
\end{picture}}

\newcommand{\tdquatreun}[4]{\begin{picture}(20,12)(-5,-1)
\put(3,0){\circle*{2}}
\put(-0.6,0){$\vee$}
\put(6,7){\circle*{2}}
\put(0,7){\circle*{2}}
\put(3,7){\circle*{2}}
\put(3,0){\line(0,1){7}}
\put(5,-2){\tiny #1}
\put(8.5,5){\tiny #2}
\put(1,10){\tiny #3}
\put(-5,5){\tiny #4}
\end{picture}}
\newcommand{\tdquatredeux}[4]{\begin{picture}(20,20)(-5,-1)
\put(3,0){\circle*{2}}
\put(-.65,0){$\vee$}
\put(6,7){\circle*{2}}
\put(0,7){\circle*{2}}
\put(0,14){\circle*{2}}
\put(0,7){\line(0,1){7}}
\put(5,-2){\tiny #1}
\put(9,5){\tiny #2}
\put(-5,5){\tiny #3}
\put(-5,12){\tiny #4}
\end{picture}}
\newcommand{\tdquatretrois}[4]{\begin{picture}(20,20)(-5,-1)
\put(3,0){\circle*{2}}
\put(-.65,0){$\vee$}
\put(6,7){\circle*{2}}
\put(0,7){\circle*{2}}
\put(6,14){\circle*{2}}
\put(6,7){\line(0,1){7}}
\put(5,-2){\tiny #1}
\put(9,5){\tiny #2}
\put(-5,5){\tiny #4}
\put(9,12){\tiny #3}
\end{picture}}
\newcommand{\tdquatrequatre}[4]{\begin{picture}(20,14)(-5,-1)
\put(3,5){\circle*{2}}
\put(-.65,5){$\vee$}
\put(6,12){\circle*{2}}
\put(0,12){\circle*{2}}
\put(3,0){\circle*{2}}
\put(3,0){\line(0,1){5}}
\put(6,-3){\tiny #1}
\put(6,4){\tiny #2}
\put(9,12){\tiny #3}
\put(-5,12){\tiny #4}
\end{picture}}
\newcommand{\tdquatrecinq}[4]{\begin{picture}(12,19)(-2,-1)
\put(0,0){\circle*{2}}
\put(0,0){\line(0,1){5}}
\put(0,5){\circle*{2}}
\put(0,5){\line(0,1){5}}
\put(0,10){\circle*{2}}
\put(0,10){\line(0,1){5}}
\put(0,15){\circle*{2}}
\put(3,-2){\tiny #1}
\put(3,3){\tiny #2}
\put(3,9){\tiny #3}
\put(3,14){\tiny #4}
\end{picture}}

\input{xy}
\xyoption{all}

\newcommand{\h}{{\cal H}}
\newcommand{\g}{\mathfrak{g}}
\newcommand{\T}{{\cal T}}

\newcommand{\D}{{\cal D}}
\newcommand{\PL}{{\cal PL}}
\newcommand{\Br}{{\cal B}r}
\newcommand{\Perm}{{\cal NAP}erm}

\title{Free brace algebras are free pre-Lie algebras}
\date{}
\author{Loïc Foissy \\
\\
{\small{\it Laboratoire de Mathématiques, FRE3111, Université de Reims}}\\
\small{{\it Moulin de la Housse - BP 1039 - 51687 REIMS Cedex 2, France}}\\
\small{e-mail : loic.foissy@univ-reims.fr}}

\newtheorem{defi}{\indent Definition}
\newtheorem{lemma}[defi]{\indent Lemma}
\newtheorem{cor}[defi]{\indent Corollary}
\newtheorem{theo}[defi]{\indent Theorem}
\newtheorem{prop}[defi]{\indent Proposition}

\newenvironment{proof}{{\bf Proof.}}{\hfill $\Box$}

\begin{document}

\maketitle

ABSTRACT. Let $\g$ be a free brace algebra. This structure implies that $\g$ is also a pre-Lie algebra and a Lie algebra.
It is already known that $\g$ is a free Lie algebra. We prove here that $\g$ is also a free pre-Lie algebra,
using a description of $\g$ with the help of planar rooted trees, a permutative product, and manipulations on the Poincaré-Hilbert series of $\g$.\\


AMS CLASSIFICATION. 17A30, 05C05, 16W30.\\

\tableofcontents

\section*{Introduction}

Let $\D$ be a set. The Connes-Kreimer Hopf algebra of rooted trees $\h_R^\D$ is introduced in \cite{Connes3}
in the context of Quantum Field Theory and Renormalization. It is a graded, connected, commutative, non-cocommutative Hopf algebra. 
If the characteristic of the base field is zero, the Cartier-Quillen-Milnor-Moore theorem insures that its dual $(\h_R^\D)^*$ is the enveloping algebra 
of a Lie algebra, based on rooted trees (note that $(\h_R^\D)^*$ is isomorphic to the Grossman-Larson Hopf algebra \cite{Grossman3,Grossman2},
as proved in \cite{Hoffman,Panaite}). This Lie algebra admits an operadic interpretation: it is the free pre-Lie algebra $\PL(\D)$ generated by $\D$, 
as shown in \cite{Chapoton}; recall that a (left) pre-Lie algebra, also called a Vinberg algebra or a left-symmetric algebra, is a vector space
$V$ with a product $\circ$ satisfying:
$$(x \circ y)\circ z-x\circ (y\circ z)=(y\circ x)\circ z-y\circ (x \circ z).$$

A non-commutative version of these objects is introduced in \cite{Foissy3,Holtkamp}. Replacing rooted trees by planar rooted trees, 
a Hopf algebra $\h_{PR}^\D$ is constructed. This self-dual Hopf algebra is isomorphic to the Loday-Ronco free dendriform algebra based 
on planar binary trees \cite{Loday}, so  by the dendriform Milnor-Moore theorem \cite{Chapoton2,Ronco}, the space of its primitive elements, 
or equivalently the space of the primitive elements of its dual, admits a structure of brace algebra, described in terms of trees in \cite{Foissy6} 
by graftings of planar forests on planar trees, and is in fact the free brace algebra $\Br(\D)$ generated by $\D$. 
This structure implies also a structure of pre-Lie algebra on $\Br(\D)$.

As a summary, the brace structure of $\Br(\D)$ implies a pre-Lie structure on $\Br(\D)$, which implies a Lie structure on $\Br(\D)$.
It is already proved in several ways that $\PL(\D)$ and $\Br(\D)$ are free Lie algebras in characteristic zero \cite{Chapoton3,Foissy6}.
A remaining question was the structure of $\Br(\D)$ as a pre-Lie algebra. The aim of the present text is to prove that $\Br(\D)$ is a free pre-Lie algebra.
We use for this the notion of non-associative permutative algebra \cite{Livernet} and a manipulation of formal series.
More precisely, we introduce in the second section of this text a non-associative permutative product $\star$ on $\Br(\D)$ and we show
that $(\Br(\D),\star)$ is free. As a corollary, we prove that the abelianisation of $\h_{PR}^\D$ (which is not $\h_R^\D$), is isomorphic to 
a Hopf algebra $\h_R^{\D'}$ for a good choice of $\D'$. This implies that $(\h_{PR}^\D)_{ab}$ is a cofree coalgebra and we recover in a different way
the result of freeness of $\Br(\D)$ as a Lie algebra in characteristic zero.
Note that a similar result for algebras with two compatible associative products is proved with the same pattern in \cite{Dotsenko}.
\\

{\bf Notations.} We denote by $K$ a commutative field of characteristic zero. All objects (vector spaces, algebras\ldots) will be taken over $K$.

\section{A description of free pre-Lie and brace algebras}

\subsection{Rooted trees and planar rooted trees}

\begin{defi} \textnormal{ \begin{enumerate}
\item A {\it rooted tree} $t$ is a finite graph, without loops, with a special vertex called the {\it root} of $t$. 
The {\it weight} of $t$ is the number of its vertices.  The set of rooted trees will be denoted by $\T$.
\item A {\it planar rooted tree} $t$  is a rooted tree with an imbedding in the plane. 
the set of planar rooted trees will be denoted by $\T_P$.
\item Let $\D$ be a nonempty set. A rooted tree decorated by $\D$ is a rooted tree
with an application from the set of its vertices into $\D$. 
The set of rooted trees decorated by $\D$ will be denoted by $\T^\D$.
\item Let $\D$ be a nonempty set. A planar rooted tree decorated by $\D$ is a planar tree
with an application from the set of its vertices into $\D$. 
The set of planar rooted trees decorated by $\D$ will be denoted by $\T_P^\D$.
\end{enumerate}} \end{defi}

{\bf Examples.} \begin{enumerate}
\item Rooted trees with weight smaller than $5$:
$$\tun,\tdeux,\ttroisun,\ttroisdeux,\tquatreun, \tquatredeux,\tquatrequatre,\tquatrecinq,
\tcinqun, \tcinqdeux,\tcinqcinq,\tcinqsept,
\tcinqneuf,\tcinqdix,\tcinqonze,\tcinqtreize,\tcinqquatorze.$$
\item Rooted trees decorated by $\D$ with weight smaller than $4$:
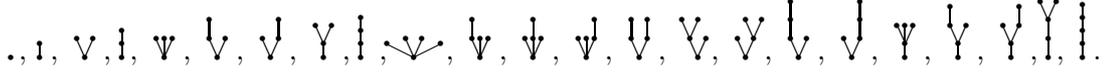
$$\tdun{$a$},\: a\in \D,\hspace{1cm} \tddeux{$a$}{$b$},\: (a,b)\in \D^2,\hspace{1cm}
\tdtroisun{$a$}{$c$}{$b$}=\tdtroisun{$a$}{$b$}{$c$},\: \tdtroisdeux{$a$}{$b$}{$c$},\:(a,b,c)\in \D^3,$$
$$ \tdquatreun{$a$}{$d$}{$c$}{$b$}=\tdquatreun{$a$}{$c$}{$d$}{$b$}
=\tdquatreun{$a$}{$d$}{$b$}{$c$}=\tdquatreun{$a$}{$b$}{$d$}{$c$}
=\tdquatreun{$a$}{$c$}{$b$}{$d$}=\tdquatreun{$a$}{$b$}{$c$}{$d$},\: \tdquatredeux{$a$}{$d$}{$b$}{$c$},\:
 \tdquatrequatre{$a$}{$b$}{$d$}{$c$}= \tdquatrequatre{$a$}{$b$}{$c$}{$d$},\:
\tdquatrecinq{$a$}{$b$}{$c$}{$d$},\:(a,b,c,d)\in \D^4.$$
\item Planar rooted trees with weight smaller than $5$:
$$\tun,\tdeux,\ttroisun,\ttroisdeux,\tquatreun, \tquatredeux,\tquatretrois,\tquatrequatre,\tquatrecinq,
\tcinqun, \tcinqdeux,\tcinqtrois,\tcinqquatre,\tcinqcinq,\tcinqsix,\tcinqsept,\tcinqhuit,
\tcinqneuf,\tcinqdix,\tcinqonze,\tcinqdouze,\tcinqtreize,\tcinqquatorze.$$
\item Planar rooted trees decorated by $\D$ with weight smaller than $4$:
$$\tdun{$a$},\: a\in \D,\hspace{1cm} \tddeux{$a$}{$b$},\: (a,b)\in \D^2,\hspace{1cm}
\tdtroisun{$a$}{$c$}{$b$},\: \tdtroisdeux{$a$}{$b$}{$c$},\:(a,b,c)\in \D^3,$$
$$ \tdquatreun{$a$}{$d$}{$c$}{$b$},\: \tdquatredeux{$a$}{$d$}{$b$}{$c$},\:
\tdquatretrois{$a$}{$c$}{$d$}{$b$},\: \tdquatrequatre{$a$}{$b$}{$d$}{$c$},\:
\tdquatrecinq{$a$}{$b$}{$c$}{$d$},\:(a,b,c,d)\in \D^4.$$
\end{enumerate}

Let $t_1,\ldots,t_n$ be elements of $\T^\D$ and let $d\in \D$. We denote by $B_d(t_1\ldots t_n)$ the rooted tree obtained by grafting 
$t_1,\ldots,t_n$ on a common root decorated by $d$. For example, $B_d(\tddeux{$a$}{$b$}\tdun{$c$})=\tdquatredeux{$d$}{$c$}{$a$}{$b$}$.
This application $B_d$ can be extended in an operator:
$$B_d: \left\{ \begin{array}{rcl}
K[\T^\D] &\longrightarrow & K\T^\D\\
t_1\ldots t_n&\longrightarrow &B_d(t_1\ldots t_n),
\end{array}\right.$$
where $K[\T^\D]$ is the polynomial algebra generated by $\T^\D$ over $K$ and $K\T^\D$ is the $K$-vector space generated by $\T^\D$.
This operator is monic, and moreover $K\T^\D$ is the direct sum of the images of the $B_d$'s, $d\in \D$.

Similarly, let $t_1,\ldots,t_n$ be elements of $\T_P^\D$ and let $d\in \D$. We denote by $B_d(t_1\ldots t_n)$ the planar rooted tree obtained by
grafting $t_1,\ldots,t_n$ in this order from left to right on a common root decorated by $d$. For example, 
$B_a(\tddeux{$b$}{$c$}\tdun{$d$})=\tdquatredeux{$a$}{$d$}{$b$}{$c$}$ 
and $B_a(\tdun{$d$}\tddeux{$b$}{$c$})=\tdquatretrois{$a$}{$b$}{$c$}{$d$}$. 
This application $B_d$ can be extended in an operator:
$$B_d: \left\{ \begin{array}{rcl}
K\langle \T_P^\D\rangle &\longrightarrow & K\T_P^\D\\
t_1\ldots t_n&\longrightarrow &B_d(t_1\ldots t_n),
\end{array}\right.$$
where $K\langle\T_P^\D\rangle$ is the free associative algebra generated by $\T_P^\D$ over $K$ and $K\T_P^\D$ is the $K$-vector space generated
by $\T_P^\D$. This operator is monic, and moreover $K\T_P^\D$ is the direct sum of the images of the $B_d$'s, $d\in \D$.

\subsection{Free pre-Lie algebras}

\begin{defi}\textnormal{A (left) pre-Lie algebra is a couple $(A,\circ)$ where $A$ is a vector space and $\circ :A\otimes A\longrightarrow A$ 
satisfying the following relation: for all $x,y,z \in A$, 
$$(x\circ y) \circ z-x \circ (y\circ z)=(y\circ x) \circ z-y \circ (x\circ z).$$}
\end{defi}

Let $\D$ be a set. A description of the free pre-Lie algebra $\PL(\D)$ generated by $\D$ is given in \cite{Chapoton}. 
As a vector space, it has a basis given by $\T^\D$, and its pre-Lie product is given, for all $t_1,t_2 \in \T^\D$, by:
$$t_1 \circ t_2=\sum_{\mbox{\scriptsize $s$ vertex of $t_2$}} \mbox{grafting of $t_1$ on $s$}.$$
For example:
$$\tdun{$a$} \circ \tdtroisun{$d$}{$c$}{$b$}
=\tdquatreun{$d$}{$c$}{$b$}{$a$}+\tdquatredeux{$d$}{$c$}{$b$}{$a$}+\tdquatretrois{$d$}{$c$}{$a$}{$b$}
=\tdquatreun{$d$}{$c$}{$b$}{$a$}+\tdquatredeux{$d$}{$c$}{$b$}{$a$}+\tdquatredeux{$d$}{$b$}{$c$}{$a$}.$$
In other terms, the pre-Lie product can be inductively defined by:
$$ \left\{ \begin{array}{rcl}
t \circ \tdun{$d$}&\longrightarrow & B_d(t),\\
t\circ B_d(t_1\ldots t_n)&\longrightarrow &\displaystyle B_d(tt_1\ldots t_n)+\sum_{i=1}^n B_d(t_1\ldots (t \circ t_i) \ldots t_n).
\end{array}\right.$$

\begin{lemma} \label{3}
Let $\D$ a set. We suppose that $\D$ has a gradation $(\D(n))_{n\in \mathbb{N}}$ such that, for all $n\in \mathbb{N}$, $\D(n)$ is finite set of 
cardinality denoted by $d_n$, and $\D(0)$ is empty. We denote by $F_\D(x)$ the Poincaré-Hilbert series of this set:
$$F_{\D}(x)=\sum_{n=1}^\infty d_n x^n.$$
This gradation induces a gradation $(\PL(\D)(n))_{n\in \mathbb{N}}$ of $\PL(\D)$. Moreover, for all $n \geq 0$,
$\PL(\D)(n)$ is finite-dimensional. We denote by $t^\D_n$ its dimension. Then the Poincaré-Hilbert series of $\PL(\D)$ satisfies:
$$F_{\PL(\D)}(x)=\sum_{n=1}^\infty t^\D_n x^n=\frac{F_{\D}(x)}{\displaystyle \prod_{i=1}^\infty (1-x^i)^{t^\D_i}}.$$
\end{lemma}

\begin{proof} The formal series of the space $K[\T^\D]$ is given by:
$$F(x)=\prod_{i=1}^\infty \frac{1}{(1-x^i)^{t^\D_i}}.$$
Moreover, for all $d\in \D(n)$, $B_d$ is homogeneous of degree $n$, so the Poincaré-Hilbert series of $Im(B_d)$
is $x^n F(x)$. As $\PL(\D)=K\T^\D=\bigoplus Im(B_d)$ as a graded vector space, its Poincaré-Hilbert formal series is:
$$F_{\PL(\D)}(x)=F(x)\sum_{n=1}^\infty d_nx^n=F(x)F_{\D}(x),$$
which gives the announced result. \end{proof}

\subsection{Free brace algebras}

\begin{defi}\textnormal{\cite{Aguiar,Chapoton2,Ronco}
A brace algebra is a couple $(A,\langle\rangle)$ where $A$ is a vector space and $\langle \rangle$ is a family of operators
$A^{\otimes n}\longrightarrow A$ defined for all $n \geq 2$:
$$\left\{ \begin{array}{rcl}
A^{\otimes n} &\longrightarrow & A\\
a_1\otimes \ldots \otimes a_n &\longrightarrow &\langle a_1,\ldots,a_{n-1};a_n\rangle,
\end{array}\right.$$
 with the following compatibilities: for all $a_1,\ldots,a_m$, $b_1,\ldots,b_n$, $c \in A$,
$$\langle a_1,\ldots,a_m;\langle b_1,\ldots, b_n; c \rangle \rangle
=\sum \langle \langle A_0,\langle A_1;b_1\rangle,A_2,\langle A_3;b_2\rangle,A_4,\ldots,A_{2n-2},\langle A_{2n-1};b_n\rangle,A_{2n};c\rangle,$$
where this sum runs over partitions of the ordered set $\{a_1,\ldots,a_n\}$ into (possibly empty) consecutive intervals
$A_0\sqcup \ldots \sqcup A_{2n}$. We use the convention $\langle a \rangle=a$ for all $a\in A$.} \end{defi}

For example, if $A$ is a brace algebra and $a,b,c \in A$:
$$\langle a; \langle b;c \rangle \rangle=\langle a,b;c \rangle+\langle b,a;c \rangle +\langle \langle a;b \rangle;c \rangle.$$
As an immediate corollary, $(A,\langle-;-\rangle)$ is a pre-Lie algebra.
Here is another example of relation in a brace algebra: for all $a,b,c,d \in A$,
$$\langle a,b;\langle c;d \rangle \rangle=\langle a,b,c;d\rangle+\langle a,\langle b;c\rangle; d\rangle
+\langle \langle a,b;c \rangle;d\rangle +\langle a,c,b;d \rangle+\langle \langle a;c \rangle, b;d\rangle
+\langle c,a,b;d\rangle.$$

Let $\D$ be a set. A description of the free brace algebra $\Br(\D)$ generated by $\D$ is given in \cite{Chapoton2,Foissy3}.
As a vector space, it has a basis given by $\T_P^\D$ and the brace structure is given, for all $t_1,\ldots,t_n \in \T_P^\D$, by:
$$\langle t_1,\ldots;t_n \rangle =\sum \mbox{graftings of $t_1\ldots t_{n-1}$ over $t_n$.}$$
Note that for any vertex $s$ of $t_n$, there are several ways of grafting a planar tree on $s$. For example:
$$\langle \tdun{$a$},\tdun{$b$};\tddeux{$d$}{$c$} \rangle
=\tdquatreun{$d$}{$c$}{$b$}{$a$}+\tdquatretrois{$d$}{$c$}{$b$}{$a$}+\tdquatreun{$d$}{$b$}{$c$}{$a$}
+\tdquatrequatre{$d$}{$c$}{$b$}{$a$}+\tdquatredeux{$d$}{$b$}{$c$}{$a$}+\tdquatreun{$d$}{$b$}{$a$}{$c$}.$$

As a consequence, the pre-Lie product of $\Br(\D)$ can be inductively defined in this way:
$$ \left\{ \begin{array}{rcl}
\langle t;\tdun{$d$} \rangle &\longrightarrow & B_d(t),\\
\langle t;B_d(t_1\ldots t_n) \rangle &\longrightarrow &\displaystyle \sum_{i=0}^n B_d(t_1\ldots t_i t t_{i+1} \ldots t_n)
+\sum_{i=1}^n B_d (t_1\ldots t_{i-1}\langle t;t_i \rangle t_{i+1}\ldots t_n).
\end{array}\right.$$

\begin{prop}
$\Br(\D)$ is the free brace algebra generated by $\D$.
\end{prop}

\begin{proof} From \cite{Chapoton2,Foissy3}. \end{proof}

\begin{lemma} \label{5}
Let $\D$ a set, with the hypotheses and notations of lemma \ref{3}. The gradation of $\D$ induces a gradation $(\Br(\D)(n))_{n\in \mathbb{N}}$ 
of $\Br(\D)$. Moreover, for all $n \geq 0$, $\Br(\D)(n)$ is finite-dimensional. Then the Poincaré-Hilbert series of $\Br(\D)$ is:
$$F_{\Br(\D)}(x)=\sum_{n=1}^\infty t'^\D_n x^n=\frac{\displaystyle 1-\sqrt{1-4F_{\D}(x)}}{2}.$$
\end{lemma}

\begin{proof} The Poincaré-Hilbert formal series of $K\langle \T_P^\D\rangle$ is given by:
$$F(x)=\frac{1}{1-F_{\Br(\D)}(x)}.$$
Moreover, for all $d\in \D(n)$, $B_d$ is homogeneous of degree $n$, so the Poincaré-Hilbert series of $Im(B_d)$
is $x^n F(x)$. As $\Br(\D)=K\T_P^\D=\bigoplus Im(B_d)$ as a graded vector space, its Poincaré-Hilbert formal series is:
$$F_{\Br(\D)}(x)=F(x) \sum_{n=1}^\infty d_nx^n=F(x)F_\D(x).$$
As a consequence, $F_{\Br(\D)}(x)-F_{\Br(\D)}(x)^2=F_\D(x),$ which implies the announced result. \end{proof}

\section{A non-associative permutative product on $\Br(\D)$}

\subsection{Definition and recalls}

The following definition is introduced in \cite{Livernet}:

\begin{defi} \textnormal{
A (left) non-associative permutative algebra is a couple $(A,\star)$, where $A$ is a vector space and $\star:A\otimes A\longrightarrow A$
satisfies the following property: for all $x,y,z \in A$,
$$x \star (y \star z)=y \star (x \star z).$$
}\end{defi}

Let $\D$ be a set. A description of the free non-associative permutative algebra $\Perm(\D)$ generated by $\D$ is given in \cite{Livernet}. As a vector space,
$\Perm(\D)$ is equal to $K\T^\D$. The non-associative permutative product is given in this way: for all $t_1 \in \T^D$, $t_2=B_d(F_2) \in \T^D$,
$$t_1 \star t_2=B_d(t_1F_2).$$
 In other terms, $t_1 \star t_2$ is the tree obtained by grafting $t_1$ on the root of $t_2$.
 As $\Perm(\D)=\PL(\D)$ as a vector space, lemma \ref{3} is still true when one replaces $\PL(\D)$ by $\Perm(\D)$.

\subsection{Permutative structures on planar rooted trees}

Let us fix now a non-empty set $\D$. We define the following product on $\Br(\D)=K\T_P^\D$:
for all $t \in \T_P^\D$, $t'=B_d(t_1\ldots t_n) \in \T_P^\D$,
$$t \star t'=\sum_{i=0}^n B_d (t_1\ldots t_i t t_{i+1}\ldots t_n).$$

\begin{prop}
$(\Br(\D),\star)$ is a non-associative permutative algebra.
\end{prop}

\begin{proof} Let us give $K\langle \T_P^\D \rangle$ its shuffle product: for all $t_1,\ldots,t_{m+n} \in \T_P^\D$,
$$(t_1\ldots t_m)*(t_{m+1}\ldots t_{m+n})=\sum_{\sigma \in Sh(m,n)} t_{\sigma^{-1}(1)} \ldots t_{\sigma^{-1}(m+n)},$$
where $Sh(m,n)$ is the set of permutations of $S_{m+n}$ which are increasing on $\{1,\ldots,m\}$ and $\{m+1,\ldots,m+n\}$.
It is well known that $*$ is an associative, commutative product. For example, for all $t,t_1,\ldots,t_n \in \T_P^\D$:
$$t *(t_1\ldots t_n)=\sum_{i=0}^n t_1\ldots t_i t t_{i+1}\ldots t_n.$$
As a consequence, for all $x \in K\T_P^\D$, $y \in K\langle \T_P^\D \rangle$, $d \in \D$:
\begin{equation}
\label{E1} x \star B_d(y)=B_d(x *y).
\end{equation}

Let $t_1,t_2,t_3=B_d(F_3) \in \T_P^\D$. Then, using (\ref{E1}):
\begin{eqnarray*}
t_1 \star (t_2 \star t_3)&=&t_1 \star B_d(t_2* F_3)\\
&=&B_d(t_1*(t_2*F_3))\\
&=&B_d((t_1*t_2)*F_3)\\
&=&B_d((t_2*t_1)*F_3)\\
&=&B_d(t_2*(t_1*F_3))\\
&=&t_2 \star (t_1 \star t_3).
\end{eqnarray*}
So $\star$ is a non-associative permutative product on $\Br(\D)$. \end{proof}

\subsection{Freeness of $\Br(\D)$ as a non-associative permutative algebra}

We now assume that $\D$ is finite, of cardinality $D$. We can then assume that $\D=\{1,\ldots,D\}$.

\begin{theo} \label{8}
$(\Br(\D),\star)$ is a free non-associative permutative algebra.
\end{theo}

\begin{proof} We graduate $\D$ by putting $\D(1)=\D$. Then $\Br(\D)$ is graded, the degree of a tree $t \in \T_P^\D$ being the number 
of its vertices. By lemma \ref{5}, as the Poincaré-Hilbert series of $\D$ is $F_\D(x)=Dx$, the Poincaré-Hilbert series of $\Br(\D)$ is:
\begin{equation}
\label {E2} F_{\Br(\D)}(x)=\sum_{i=1}^\infty t'^\D_i x^i=\frac{1-\sqrt{1-4Dx}}{2}.
\end{equation}
We consider the following isomorphism of vector spaces:
$$B: \left\{ \begin{array}{rcl}
(K\langle \T_P^\D\rangle)^d &\longrightarrow & \Br(\D)\\
(F_1,\ldots,F_D)&\longrightarrow & \displaystyle \sum_{i=1}^d B_i(F_i).
\end{array}\right.$$
Let us fix a graded complement $V$ of the graded subspace $\Br(\D) \star \Br(\D)$ in $\Br(\D)$.
Because $\Br(\D)$ is a graded and connected (that is to say $\Br(\D)(0)=(0)$), $V$ generates $\Br(\D)$ as a non-associative permutative algebra.
By (\ref{E1}), $\Br(\D)\star \Br(\D)=B((\T_P^\D * K\langle \T_P^\D \rangle)^D)$.

Let us then consider $\T_P^\D * K\langle \T_P^\D\rangle$, that is to say the ideal of $(K\langle \T_P^\D \rangle,*)$ generated by $\T_P^\D$.
It is known that $(K\langle \T_P^\D \rangle,*)$ is isomorphic to a symmetric algebra (see \cite{Ree}). Hence, there exists a graded subspace $W$
of $K\langle \T_P^\D \rangle$, such that $(K\langle \T_P^\D \rangle,*)\approx S(W)$ as a graded algebra. We can assume that $W$ contains
$K\T_P^\D$. As a consequence:
\begin{equation}
\label{E3} \frac{K\langle \T_P^\D \rangle}{\T_P^\D*K\langle \T_P^\D \rangle}
\approx \frac{S(W)}{S(W)\T_P^\D}\approx S\left(\frac{W}{K\T_P^\D} \right).
\end{equation}
We denote by $w_i$ the dimension of $W(i)$ for all $i\in \mathbb{N}$. 
Then, the Poincaré-Hilbert formal series of $S\left(\frac{W}{K\T_P^\D} \right)$ is:
\begin{equation}
\label{E4} F_{S\left(\frac{W}{K\T_P^\D} \right)}(x)=\prod_{i=1}^\infty \frac{1}{(1-x^i)^{w_i-t'^\D_i}}.
\end{equation}
Moreover, the  Poincaré-Hilbert formal series of $K\langle \T_P^\D \rangle\approx S(W)$ is, by (\ref{E2}):
\begin{equation}
\label{E5} F_{S(W)}(x)=\frac{1}{1-F_{\Br(\D)}(x)}=\frac{1-\sqrt{1-4Dx}}{2Dx}=\frac{F_{\Br(\D)}(x)}{Dx}=\prod_{i=1}^\infty \frac{1}{(1-x^i)^{w_i}}.
\end{equation}
So, from (\ref{E3}), using (\ref{E4}) and (\ref{E5}), the Poincaré-Hilbert series of $\T_P^\D*K\langle \T_P^\D \rangle$ is:
\begin{eqnarray*}
F_{\T_P^\D*K\langle \T_P^\D \rangle}(x)&=&F_{S(W)}(x)-F_{S\left(\frac{W}{K\T_P^\D} \right)}(x)\\
&=&\prod_{i=1}^\infty \frac{1}{(1-x^i)^{w_i}}\left( 1-\prod_{i=1}^\infty (1-x^i)^{t'^\D_i}\right)\\
&=&\frac{F_{\Br(\D)}(x)}{Dx}\left( 1-\prod_{i=1}^\infty (1-x^i)^{t'^\D_i}\right).
\end{eqnarray*}
As $B$ is homogeneous of degree $1$, the Poincaré-Hilbert formal series of $\Br(\D) \star \Br(\D)$ is:
$$F_{\Br(\D)\star \Br(\D)}(x)=Dx F_{\T_P^\D*K\langle \T_P^\D \rangle}(x)=F_{\Br(\D)}(x)\left( 1-\prod_{i=1}^\infty (1-x^i)^{t'^\D_i}\right).$$
Finally, the Poincaré-Hilbert formal series of $V$ is:
$$\ F_V(x)=F_{\Br(\D)}(x)-F_{\Br(\D)\star \Br(\D)}(x)=F_{\Br(\D)}(x) \prod_{i=1}^\infty (1-x^i)^{t'^\D_i}.$$

Let us now fix a basis $(v_i)_{i\in I}$ of $V$, formed of homogeneous elements. There is a unique epimorphism of non-associative permutative algebras:
$$\Theta: \left\{ \begin{array}{rcl}
\Perm(I)&\longrightarrow & \Br(\D)\\
\tdun{$i$} &\longrightarrow & v_i.
\end{array}\right.$$
We give to $i\in I$ the degree of $v_i\in \Br(\D)$. With the induced gradation of $\Perm(I)$, $\Theta$ is a graded epimorphism.
In order to prove that it is an isomorphism, it is enough to prove that the Poincaré-Hilbert series of $\Perm(I)$ and $\Br(\D)$ are equal.
By lemma \ref{3}, the formal series of $\Perm(I)$, or, equivalently, of $\PL(I)$, is:
\begin{equation}
\label{E6} F_{\Perm(I)}(x)=\sum_{n=1}^\infty t^\D_ix^i=\frac{F_V(x)}{\displaystyle \prod_{i=1}^\infty (1-x^i)^{t^\D_i}}
=F_{\Br(\D)}(x) \prod_{i=1}^\infty (1-x^i)^{t'^\D_i-t^\D_i}.
\end{equation}
Let us prove inductively that $t_n=t'_n$ for all $n \in \mathbb{N}$. It is immediate if $n=0$, as $t_0=t'_0=0$.
Let us assume that $t^\D_i=t'^\D_i$ for all $i<n$. Then:
$$ \prod_{i=1}^\infty (1-x^i)^{t^\D_i-t'^\D_i}=1+\mathcal{O}(x^n).$$
As $t'_0=0$, the coefficient of $x^n$ in (\ref{E6}) is $t_n=t'_n$. So $F_{\Perm(I)}(x)=F_{S(W)}(x)$, and $\Theta$ is an isomorphism. \end{proof}

\section{Freeness of $\Br(\D)$ as a pre-Lie algebra}

\subsection{Main theorem}

\begin{theo}
Let $\D$ be a finite set. Then $\Br(\D)$ is a free pre-Lie algebra.
\end{theo}

\begin{proof} We give a $\mathbb{N}^2$-gradation on $\Br(\D)$ in the following way:
$$\Br(\D)(k,l)=Vect(t\in \T_P^\D\:/\: \mbox{$t$ has $k$ vertices and the fertility of its root is $l$}).$$
The following points are easy:
\begin{enumerate}
\item For all $i,j,k,l \in \mathbb{N}$, $\Br(\D)(i,j)\star \Br(\D)(k,l)\subseteq \Br(\D)(i+k,l+1)$.
\item For all $i,j,k,l \in \mathbb{N}$, $t_1 \in \Br(\D)(i,j)$, $t_2 \in \Br(\D)(k,l)$, $\langle t_1;t_2 \rangle-t_1 \star t_2 \in Br(\D)(i+k,l)$.
\end{enumerate}
Let us fix a complement $V$ of $\Br(\D) \star \Br(\D)$ in $\Br(\D)$ which is $\mathbb{N}^2$-graded. Then $\Br(\D)$ is isomorphic 
as a $\mathbb{N}$-graded non-associative permutative algebra to $\Perm(V)$, the free non-associative permutative algebra generated by $V$.

Let us prove that $V$ also generates $\Br(\D)$ as a pre-Lie algebra. As $\Br(\D)$ is $\mathbb{N}$-graded, with $\Br(\D)(0)$, it is enough to prove
that $\Br(\D)=V+\langle \Br(\D);\Br(\D) \rangle$. Let $x\in \Br(\D)(k,l)$, let us show that $x\in V+\langle \Br(\D);\Br(\D) \rangle$ by induction on $l$.
If $l=0$, then $t \in \Br(\D)(1)=V(1)$. If $l=1$,  we can suppose that $x=B_d(t)$, where $t\in \T_P^\D$. 
Then $x=\langle t;\tdun{$d$}\rangle \in \langle \Br(\D);\Br(\D) \rangle$.  Let us assume the result for all $l'<l$. 
As $V$ generates $(\Br(\D),\star)$, we can write $x$ as:
 $$x=x'+\sum_i  x_i \star y_i,$$
 where $x' \in V$ and $x_i,y_i \in \Br(\D)$. By the first point, we can assume that:
$$\sum_i x_i \otimes y_i \in \bigoplus_{i+j=k} \Br(\D)(i) \otimes \Br(\D)(j,l-1).$$
So, by the second point:
\begin{eqnarray*}
x-x'-\sum_i \langle x_i;y_i \rangle&=&\sum_i x_i \star y_i -\langle x_i;y_i \rangle\\
&\in& \sum_{i+j=k} \Br(\D)(i+j,l-1)\\
&\in& V+\langle \Br(\D);\Br(\D) \rangle,
\end{eqnarray*}
by the induction hypothesis. So $x \in V+\langle \Br(\D);\Br(\D) \rangle$.
 
 Hence, there is an homogeneous epimorphism:
 $$ \left\{ \begin{array}{rcl}
\PL(V)&\longrightarrow & \Br(\D)\\
v \in V&\longrightarrow & v.
\end{array}\right.$$
As $\PL(V)$, $\Perm(V)$ and $\Br(\D)$ have the same Poincaré-Hilbert formal series, this is an isomorphism. \end{proof}\\

We now give the number of generators of $\Br(\D)$ in degree $n$ when $card(\D)=D$ for small values of $n$, computed using lemmas \ref{3} 
and \ref{5}:
\begin{enumerate}
\item For $n=1$, $D$.
\item For $n=2$, $0$.
\item For $n=3$, $\displaystyle \frac{D^2(D-1)}{2}$.
\item For $n=4$, $\displaystyle \frac{D^2(2D-1)(2D+1)}{3}$.
\item For $n=5$, $\displaystyle \frac{D^2(31D^3-2D^2-3D-2)}{8}$.
\item For $n=6$, $\displaystyle \frac{D^2(356D^4-20D^3-5D^2+5D-6)}{30}$.
\item For $n=7$, $\displaystyle \frac{D^2(5441D^5-279D^4-91D^3-129D^2-22D-24)}{144}$.
\end{enumerate}

\subsection{Corollaries}

\begin{cor}
Let $\D$ be any set. Then $\Br(\D)$ is a free pre-Lie algebra.
\end{cor}

\begin{proof} We graduate $\Br(\D)$ by putting all the $\tdun{$d$}$'s homogeneous of degree $1$.
Let $V$ be a graded complement of $\langle \Br(\D),\Br(\D)\rangle$. There exists an epimorphism of graded pre-Lie algebras:
$$ \Theta: \left\{ \begin{array}{rcl}
\PL(V)&\longrightarrow & \Br(\D)\\
\tdun{$v$}&\longrightarrow & v.
\end{array}\right.$$
Let $x$ be in the kernel of $\Theta$. There exists a finite subset $\D'$ of $\D$, such that all the decorations of the vertices of the trees appearing in
$x$ belong to $\Br(\D')$. By the preceding theorem, as $\Br(\D')$ is a free pre-Lie algebra, $x=0$. So $\Theta$ is an isomorphism. \end{proof}

\begin{cor}
Let $\D$ be a graded set, satisfying the conditions of lemma \ref{3}. There exists a graded set $\D'$,  such that $(\h_{PR}^\D)_{ab}$ is isomorphic, 
as a graded Hopf algebra, to $\h_R^{\D'}$.
\end{cor}

\begin{proof} $(\h_{PR}^\D)_{ab}$ is isomorphic, as a graded Hopf algebra, to ${\cal U}(\Br(\D))^*$. 
For a good choice of $\D'$, $\Br(\D)$ is isomorphic to $\PL(\D')$ as a pre-Lie algebra, so also as a Lie algebra.
So ${\cal U}(\Br(\D))$ is isomorphic to ${\cal U}(\PL(\D'))$. Dually, $(\h_{PR}^\D)_{ab}$ is isomorphic to $\h_R^{\D'}$. \end{proof}

\begin{cor}
Let $\D$ be graded set, satisfying the conditions of lemma \ref{3}. Then $(\h_{PR}^\D)_{ab}$ is a cofree coalgebra. Moreover, $\Br(\D)$ is free as a Lie algebra.
\end{cor}

\begin{proof} It is proved in \cite{Foissy} that $(\h_R^{\D'})^*$ is a free algebra, so $Prim((\h_R^{\D'})^*)=\PL(\D')$ is a free Lie algebra and
$\h_R^{\D'}$ is a cofree coalgebra. So $Prim((\h_{PR}^\D)^*)=\Br(\D)$ is a free Lie algebra and $\h_{PR}^\D$ is a cofree coalgebra. \end{proof}

\bibliographystyle{amsplain}
\bibliography{biblio}

\providecommand{\bysame}{\leavevmode\hbox to3em{\hrulefill}\thinspace}
\providecommand{\MR}{\relax\ifhmode\unskip\space\fi MR }
\providecommand{\MRhref}[2]{%
  \href{http://www.ams.org/mathscinet-getitem?mr=#1}{#2}
}
\providecommand{\href}[2]{#2}
\begin{thebibliography}{10}

\bibitem{Aguiar}
Marcelo Aguiar, \emph{Infinitesimal bialgebras, pre-{L}ie and dendriform
  algebras}, Hopf algebras, Lecture Notes in Pure and Appl. Math., vol. 237,
  Dekker, 2004, arXiv:math/0211074, pp.~1--33.

\bibitem{Chapoton2}
Frédéric Chapoton, \emph{Un théorème de {C}artier-{M}ilnor-{M}oore-{Q}uillen
  pour les bigèbres dendriformes et les algèbres braces}, J. Pure Appl. Algebra
  \textbf{168} (2002), no.~1, 1--18, arXiv:math/0005253.

\bibitem{Chapoton3}
\bysame, \emph{Free pre-lie algebras are free as lie algebras},
  arXiv:0704.2153, to appear in {\it Bulletin canadien de mathématiques}, 2007.

\bibitem{Chapoton}
Frédéric Chapoton and Muriel Livernet, \emph{Pre-{L}ie algebras and the rooted
  trees operad}, Internat. Math. Res. Notices \textbf{8} (2001), 395--408,
  arXiv:math/0002069.

\bibitem{Connes3}
Alain Connes and Dirk Kreimer, \emph{Hopf algebras, {R}enormalization and
  {N}oncommutative geometry}, Comm. Math. Phys \textbf{199} (1998), no.~1,
  203--242, arXiv:hep-th/9808042.

\bibitem{Dotsenko}
Vladimir Dotsenko, \emph{Compatible associative products and trees},
  arXiv:0809.1773, to appear in Algebra and Number Theory, 2008.

\bibitem{Foissy}
Lo{\"\i}c Foissy, \emph{Finite-dimensional comodules over the {H}opf algebra of
  rooted trees}, J. Algebra \textbf{255} (2002), no.~1, 85--120,
  arXiv:math.QA/0105210.

\bibitem{Foissy6}
\bysame, \emph{Les alg\`ebres de {H}opf des arbres enracin\'es, {I}}, Bull.
  Sci. Math. \textbf{126} (2002), no.~3, 193--239, arXiv:math.QA/0105212.

\bibitem{Foissy3}
\bysame, \emph{Les alg\`ebres de {H}opf des arbres enracin\'es, {II}}, Bull.
  Sci. Math. \textbf{126} (2002), no.~4, 1249--288, arXiv:math.QA/0105212.

\bibitem{Grossman3}
Robert~L. Grossman and Richard~G. Larson, \emph{Hopf-algebraic structure of
  families of trees}, J. Algebra \textbf{126} (1989), no.~1, 184--210,
  arXiv:0711.3877.

\bibitem{Grossman2}
\bysame, \emph{Hopf-algebraic structure of combinatorial objects and
  differential operators}, Israel J. Math. \textbf{72} (1990), no.~1-2,
  109--117.

\bibitem{Hoffman}
Michael~E. Hoffman, \emph{Combinatorics of rooted trees and {H}opf algebras},
  Trans. Amer. Math. Soc. \textbf{355} (2003), no.~9, 3795--3811.

\bibitem{Holtkamp}
Ralf Holtkamp, \emph{Comparison of {H}opf algebras on trees}, Arch. Math.
  (Basel) \textbf{80} (2003), no.~4, 368--383.

\bibitem{Livernet}
Muriel Livernet, \emph{A rigidity theorem for pre-{L}ie algebras}, J. Pure
  Appl. Algebra \textbf{207} (2006), no.~1, 1--18, arXiv:math/0504296.

\bibitem{Loday}
Jean-Louis Loday and Maria~O. Ronco, \emph{Hopf algebra of the planar binary
  trees}, Adv. Math. \textbf{139} (1998), no.~2, 293--309.

\bibitem{Panaite}
Florin Panaite, \emph{Relating the {C}onnes-{K}reimer and {G}rossman-{L}arson
  {H}opf algebras built on rooted trees}, Lett. Math. Phys. \textbf{51} (2000),
  no.~3, 211--219.

\bibitem{Ree}
Rimhak Ree, \emph{Lie elements and an algebra associated with shuffles}, Ann.
  of Math. (2) \textbf{68} (1958), 210--220.

\bibitem{Ronco}
Maria Ronco, \emph{A {M}ilnor-{M}oore theorem for dendriform {H}opf algebras},
  C. R. Acad. Sci. Paris Sér. I Math. \textbf{332} (2001), no.~2, 109--114.

\end{thebibliography}

\end{document}